\documentclass[12pt, a4paper]{article}

\usepackage{amsmath, amsfonts, amsthm, amssymb}

\newcommand     {\C}            {{\mathbb C}}
\newcommand     {\F}            {{\mathbb F}}

\newcommand     {\HP}           {{\mathbb H}}
\newcommand     {\Ind}          {\sym{Ind}}
\newcommand     {\Inv}          {\sym{Inv}}

\newcommand     {\Mpp}[1][\Z]   {\sym{Mp}(2,{#1})}
\newcommand     {\Q}            {{\mathbb Q}}

\newcommand     {\SL}[1][\Z]    {\sym{SL}(2,{#1})}
\newcommand     {\Z}            {{\mathbb Z}}

\newcommand     {\e}            {\sym{e}}
\newcommand     {\conj}[1]      {{#1}^c}

\newcommand     {\pt}[1]        {\sym{e}_{#1}}
\newcommand     {\sym}[1]       {\operatorname{#1}}

\theoremstyle{plain}
\newtheorem{Theorem}{Theorem}

\newtheorem*{Lemma}{Lemma}
\theoremstyle{definition}
\newtheorem*{Definition}{Definition}

\title{%
  Jacobi Forms of Degree One\\ and Weil Representations
}

\author{%
  Nils-Peter Skoruppa}

\date{}

\begin{document}

\maketitle

\begin{abstract}
\noindent
We discuss the notion of Jacobi forms of degree one
with matrix index,
we state dimension formulas, give explicit examples,
and indicate how closely their theory is connected to
the theory of invariants of Weil representations
associated to finite quadratic modules.
\end{abstract}

\section{Jacobi forms of degree one}
\label{Jacobi}

Jacobi forms of degree one with matrix index $F$ gained recent
interest, mainly due to applications in the theory of Siegel and
orthogonal modular forms, and in the
geometry of moduli spaces.  Of
particular interest among these Jacobi forms are those of {\it critical weight},
i.e.~those whose weight equals $(\sym{rank}(F)+1)/2$. There are no
Jacobi forms of of index~$F$ and weight strictly less than
$\sym{rank}(F)/2$, and for weights strictly greater than
$(\sym{rank}(F)+3)/2$ we have at least an explicit and easily
computable dimension formula (see Theorem~\ref{thm:dimension-formula}
below).

From the point of view of algebraic geometry Jacobi forms of degree
one are holomorphic functions $\phi(\tau,z)$ of a variable $\tau$ in
the Poincar\'e upper half plane $\HP$ and of complex variables $z\in\C^n$
such that, for fixed $\tau$, the function $\phi(\tau,\cdot)$ is a
theta function on the algebraic variety $\C^n/\Lambda_\tau$, where
$\Lambda_\tau$ denotes the lattice $\Z^n\tau+\Z^n$, and such that, for
any $\tau$ and all $A$ in a subgroup of finite index in $\SL$, the
theta function $\phi(\tau,\cdot)$ on $\C^n/\Lambda_\tau$ and
$\phi(A\tau,\cdot)$ on the isomorphic torus $\C^n/\Lambda_{A\tau}$ are
related by a simple transformation formula. Thus, for fixed $\tau$,
the function $\phi(\tau,\cdot)$ corresponds to a holomorphic section
of a positive line bundle on $\C^n/\Lambda_{\tau}$. The positive line
bundles on $X_\tau=\C^n/\Lambda_{\tau}$ are (up to
translation and isomorphism) parameterized by their Chern classes in
$H^2(X_\tau,\Z)$. It is not difficult to show that the cone of positive
Chern classes in $H^2(X_\tau,\Z)$ is in one to one correspondence with
the set of symmetric, positive definite matrices $F$ with entries in
$\frac12 \Z$ via the map
$$
F=\left(f_{p,q}\right)
\mapsto
\frac 1{i\Im(\tau)}\sum_{p,q}f_{p.q}\,dz_p \wedge d\overline{z}_q
,
$$
where $z_p$ denote the coordinate functions associated to the
canonical basis of~$\C^n$. This explains the appearance and the nature
of the matrix index in the formal definition of a Jacobi form, and it
shows also that the usual theory of Jacobi forms omits the case of
those indices $F$ where the diagonal entries of $F$ are not
necessarily integral (as we shall do too in the following discussion;
however, see the remarks after the formal definition of the notion
{\it Jacobi form} in the Appendix).  Recall that a symmetric, positive
definite matrix $F$ is called {\it half integral} if it has
half integral entries, but integers on the diagonal.

For the formal definition of the space $J_{k,F}(\Gamma,\chi)$ of
Jacobi forms of integral weight $k$, of positive-definite
half integral index $F$, on a subgroup $\Gamma$ of finite index in
$\SL$, and with character $\chi$ we refer the reader to the
Appendix. The implicit restrictions made at this point can be relaxed:
One can admit half integral $k$, semi-positive definite $F$ with
half-integers on the diagonal and vector-valued Jacobi-forms. We
suppress the discussion of these possible generalizations for not
overloading this presentation by too many technical details; the
interested reader is referred to \cite{Sko 07} for a more general and
elaborated treatment.

There is an explicit  dimension formula for $J_{k,F}(\Gamma,\chi)$ if
$k > (n+3)/2$, where $n$ denotes the rank of $F$. For understanding
why such a formula exists note that, for fixed $\tau$,
the dimension of the space
of holomorphic sections of the line bundle on $\C^n/\Lambda_\tau$
corresponding to Jacobi forms of index $F$ has dimension $\det(2F)$. A
basis for the space of theta functions corresponding to this line
bundle is given by the special theta functions
$$
\vartheta_{F,x}(\tau,z)
	 =
	 \sum_{\begin{subarray}{c}r\in\Z^n\\
		  r\equiv x \bmod 2F\Z^n\end{subarray}}
	 \e\big(\tau \frac14 F^{-1}[r] + r^tz\big)
	 \qquad
(x\in\Z^n)
.
$$
Note that $\vartheta_{F,x}$ depends only on $x$ modulo $2F\Z^n$. Thus
any Jacobi $\phi$ form of index $F$ can be written in the form
$$
\psi(\tau,z)=\sum_{x\in\Z^n/2F\Z^n}
h_x(\tau)\,\vartheta_{F,x}(\tau,z)
$$
with functions $h_x$ which are holomorphic in $\HP$.  The invariance
of $\phi$ under $\Gamma$ is, in this representation, reflected by the
fact that the $h_x$ are the coordinates of a vector
valued modular form associated to the (dual of the) Weil
representation $W(F)$ which we shall explain in section~\ref{Weil}. If
the weight of $\phi$ is $k$ then the weight of the corresponding
vector valued modular form is $k-\frac n2$. This shows already that
there are no nontrivial Jacobi forms of index $F$
and weight $k<\frac n2$. In any
case we can apply the Eichler-Selberg trace formula or Shimura's
variant (based on the Lefschez fixed point theorem) to obtain a
dimension formula for vector valued modular forms
\cite[Sec.~4.3]{Eh-S 95},~\cite[p.~100]{Sko 85}, and we can then  deduce
from this a general dimension formula for Jacobi forms.

\begin{Theorem}[\cite{Sko 07}]
  \label{thm:dimension-formula}
  Let $F$ be a half integral positive definite $n\times n$ matrix, let
  $k\in\Z$, let $\Gamma$ be a subgroup of $\SL$ and let $\chi$ be
  a linear character
  of~$\Gamma$ of finite order.
  Then one has
  \begin{equation*}
	 \begin{split}
		\dim &J_{k,F}(\Gamma,\chi) 
		-
		\dim M_{2+\frac n2 - k}^{\text cusp}\otimes_{\C[\Mpp]} \conj{X(i^{n-2k})}\\
		&=
		\frac{k - \frac n2 - 1}{12}\,\dim X(i^{n-2k})
		+\frac 14
		\sym{Re}\left(e^{\pi i (k-\frac n2 )/2}\sym{tr}((S,w_S),X(i^{n-2k}))\right)\\
		&+\frac 2{3\sqrt 3}
		\sym{Re}\left(e^{\pi i (2k - n + 1)/6}\sym{tr}((ST,w_{ST}),X(i^{n-2k}))\right)
		-\sum_{j=1}^r(\lambda_j-\frac12)
		.
	 \end{split}
  \end{equation*}
  Here $S=\begin{pmatrix}0&-1\\1&0\end{pmatrix}$,
  $T=\begin{pmatrix}1&1\\0&1\end{pmatrix}$ and, for any $A$ in $\SL$,
  we use $(A,w_A)$ for the corresponding element of $\Mpp$ (whose
  precise definition is given in
  the Appendix).  Moreover, $X(i^{n-2k})$ denotes the $\Mpp$-submodule
  of all $v$ in $\conj{W(F)}\otimes\Ind_\Gamma^{\SL} \C(\chi)$ such
  that $(-1,i)v=i^{n-2k}v$, and the $\lambda_j$ are rational
  numbers $0\le \lambda_j < 1$ such that $\prod_{j=1}^r(t-e^{2\pi i\lambda_j})\in\C[t]$
  equals the characteristic polynomial of the automorphism of
  $X(i^{n-2k})$ given by $v\mapsto (T,1)v$.
\end{Theorem}

(There are some more technical explanations necessary if
the reader wants to apply the theorem.
In its statement we used the following notations from representation
theory: $\C(\chi)$ is the $\Gamma$-module with underlying vector space
$\C$ and with the action given by $\left(A,z\right)\mapsto \chi(A)\cdot z$, and
$\Ind_\Gamma^{\SL} \C(\chi)$ denotes the $\SL$-module induced by
$\C(\chi)$.  Moreover, $\conj{W(F)}$ is the $\Mpp$-module with
underlying vector space equal to the dual vector space of $W(F)$, and
with the $\Mpp$-action $(A,f)\mapsto Af$, where $(Af)(v)=f(A^{-1}v)$ for
all $v$ in~$W(F)$. If the
action of $\Mpp$ on $\conj{W(F)}$ does not factor through $\SL$ then,
for forming the tensor product considered in the theorem, the induced
representation has to be considered as $\Mpp$-module; otherwise we can
form the tensor product with respect to $\C[\SL]$.
The action of $\Mpp$ on $W(F)$
factors through $\SL$ if and only if $n$ is even since then
the sigma-invariant of the determinant module $D_F$ is a fourth root
of unity by Milgram's formula; cf.~section~\ref{Weil} for an
explanation of these terms.)

If $k\ge\frac n2 + 2$ the theorem gives us a ready to compute formula
since then then the second term on the left hand side vanishes. There
remain two weights where the theorem is of no help: $k=\frac n2$,
$k=\frac n2 +1$ if $n$ is even, and $k=\frac {n+1}2$ (critical weigh),
$k=\frac{n+3}2$ if $n$ is odd. For the case $k=\frac n2 +1$ with even
$n$ we do not know of any approach to set up a general trace formula;
in fact, here the underlying vector valued modular forms are of
weight one. The case $k=\frac{n+3}2$ with odd $n$ can be treated by
methods similar to the one for critical weight, which we shall explain
in a moment.  The cases $\frac n2$ and $\frac{n+1}2$ for even and odd
weights, respectively, can be reduced to purely representation
theoretic questions concerning the Weil representation $W(F)$ (for the
critical weight case we even have to assume additionally that the
kernel of $\chi$ is a congruence subgroup).

Assume first of all that $n$ is even and $k=\frac n2$. Then the $h_x$
introduced above have weight 0 and are consequently constants. This leads to
the following theorem.

\begin{Theorem}[\cite{Sko 07}]
  \label{thm:singular-weight}
  There exists a natural isomorphism
  $$
  J_{\frac n2,F}(\Gamma,\chi)
  \cong
  W(F)^*\otimes_{\C[\SL]}\Ind_\Gamma^{\SL}\C(\chi)
  .
  $$
\end{Theorem}
(Here $ W(F)^*$ denotes the $\Mpp$-right module whose underlying vector
space is the dual vector space of $W(F)$ and the $\Mpp$-action is
given by $(f,A)\mapsto A^{-1}f$; cf the remark after
Theorem~\ref{thm:dimension-formula}.)  Note that, for trivial
$\chi$ and, say, $\Gamma$ equal to the full modular group, the right
hand side of the stated isomorphism is isomorphic to the subspace
$\Inv\left(W(F)\right)$ of $\SL$-invariants of the representation
$W(F)$. The question of describing these invariants or computing the
dimension of the space of invariants is in general unsolved. We shall
discuss this question in the second part of this article.

For odd $n$ and $k=\frac{n+1}2$, the critical weight, the description
of $J_{k,F}(\Gamma,\chi)$ becomes even more subtle. Here the $h_x$ are
modular form of weight $\frac 12$. If we assume that the kernel of
$\chi$ (and hence $\Gamma$) is a congruence subgroup the~$h_x$ are
modular forms on congruence subgroups. By a theorem of Serre-Stark the
only modular forms of weight $\frac 12$ on congruence subgroups are
theta series (more precisely, linear combinations of the null values
$\vartheta_{m,x}(\tau,0)$, where $x$ and $m$ run through the
integers and positive integers, respectively).  Based on this description, the decomposition of the
$\SL$-module of all modular forms of weight $\frac12$ (on congruence
subgroups) was explicitly derived in \cite{Sko 85}, \cite{Sko
07}. From this description we can then finally deduce the following theorem.

\begin{Theorem}[\cite{Sko 07}]
  \label{thm:second-isomorphism}
  Let $\Gamma$ be a subgroup of finite index in $\SL$, let~$\chi$ be
  a linear character of $\Gamma$, and let $F$ be half integral of size
  $n$ and level $f$.  Assume that, for some $m$, the group $\chi$ is
  trivial on $\Gamma(4m)$, and that $f$ divides $4m$. Then there is a
  natural isomorphism
  $$
	 J_{\frac {n+1}2,F}(\Gamma,\chi)
	 \longrightarrow
	 \bigoplus_{\begin{subarray}{c}l|m\\
		  m/l \text{ squarefree}\end{subarray}}
	 \Big(W\big(l\oplus F\big)^{\iota}\Big)^*\otimes_{\C[\SL]}\Ind_\Gamma^{\SL}\C(\chi)
	 $$
  Here $W\big(l\oplus F\big)^{\iota}$
  denotes the $+1$-eigenspace of the involution $\iota$ on
  $W\big(l\oplus F\big)$ induced by the automorphism
  $(x,y)\mapsto(-x,y)$ of $D_{(l)\oplus F}$.
\end{Theorem}
(Recall that the level $f$ of $F$ is the smallest positive integer
such that $f(2F)^{-1}$ is again half integral. For the notion of the
{\it determinant module} $D_G$ associated to a half integral $G$ see
section~\ref{Weil}.  By $l\oplus F$ we denote the matrix of size
$(n+1)\times (n+1)$ consisting of the two blocks $l$ and $F$ centered
on the diagonal and having
zero entries otherwise.) We observe again that, for trivial~$\chi$, the
calculation of $J_{\frac{n+1}2,F}(\Gamma)$ is basically equivalent
to the calculation
of the subspaces $\Inv\left(W(l\oplus F)\right)$.

We conclude this section by some remarks concerning the graded algebra
$$
J_{*,F}(\Gamma,\chi) = \bigoplus_{k\in\Z}J_{k,F}(\Gamma,\chi)
.
$$
Here $J_{*,F}(\Gamma,\chi)$ may be considered as the subspace of
functions on $\HP\times\C^n$ spanned by the Jacobi forms in
$J_{k,F}(\Gamma,\chi)$, where $k$ runs through all integers,
and the proof of the direct sum
decomposition is then an easy exercise. It is not hard to see from
the dimension formula of Theorem~\ref{thm:dimension-formula} that the
Hilbert-Poincar\'e series of this graded algebra is of the form
$$
\sum_{k\in\Z}
\dim J_{k,F}(\Gamma,\chi)\,x^k
=
\frac{p_F(x)}{(1-x^4)(1-x^6)}
,
$$ where $p_F(x)$ is a polynomial. Indeed, for $k\ge \frac n2+2$ the
Hilbert-Poincar\'e series splits up into four sums $S_j=\sum_k
a_j(k)x^k$ ($j=1,2,3,4$) according to the four terms on the right hand
side of Theorem~\ref{thm:dimension-formula}.  Here $a_1(k)$ is a
linear function in $k$ times a sequence which depends only on $k$
modulo $2$, whence $S_1=\text{polynomial}/(1-x^2)^2$. Moreover
$a_2(k)$, $a_3(k)$ and $a_4(k)$ depend only on~$k$ modulo $4$, $6$ and
$2$, respectively, which yields~$S_j=\text{polynomial}/(1-x^{j'})$
with $j'=4,6,2$, respectively. From this argument it is also easy to
see that the degree of $p_F(x)$ is less than or equal to $12$. We
note that $p_F(1)=\det(2F)\cdot[\SL:\Gamma]$. Namely, multiplying the
Hilbert-Poincar\'e series by $24(1-x)^2$ and letting $x$ tend to $1$
gives, on using the decomposition into the parts $S_j$ and observing
that the contributions from $S_j$ vanish for $j=2,3,4$, the value
$\dim X(i^n)+\dim X(i^n-2)=\dim \conj{W(F)}\otimes \Ind \C(\chi)$.  On
the other hand, $p_F(x)(1-x)^2/(1-x^4)(1-x^6)$ becomes $p_F(1)/24$
for~$x=1$.

The explanation for the shape of the Hilbert-Poincar\'e series of
$J_{*,F}(\Gamma,\chi)$ is as follows.  Multiplication by usual
elliptic modular forms turns $J_{*,F}(\Gamma,\chi)$ into a graded
module over the graded ring $M_*=\C[E_4,E_6]$, where $E_4, E_6$ are
the usual Eisenstein series of weight $4$ and $6$ on the full modular
group. One may copy the proof in \cite{E-Z 85} to show that, in fact,
$J_{*,F}(\Gamma,\chi)$  is free over $M_*$. If $S$ is a system of
homogeneous generators of $J_{*,F}(\Gamma,\chi)$ as module over $M_*$,
and if $s_j$ denotes the number of generators in $S$ of weight $j$
then $p_F(x)=s_0+s_1x+\cdots+s_{12}x^{12}$. In particular, the rank of
$J_{*,F}(\Gamma,\chi)$ over $M_*$ equals
$p_F(1)=\det(2F)\,[\SL:\Gamma]$.

As an example we consider $J_{k,F}:=J_{k,F}(\SL)$ for half integral $2\times2$ matrices
with $\det(2F)=p$ for primes $p\not=3$.  The dimension formula, for
$k\ge 3$, becomes then
\begin{gather*}
\dim J_{k,F}=
\frac{k-2}{12}\,\cdot\,\frac{(p+(-1)^k)}2
-\frac14 \left(\frac{-4}{k-2}\right)\left(\frac {-2}p\right)\\
-\frac13\left(\frac {k-2}3\right)\frac{\left(\frac p3\right)+(-1)^k}2
-\frac{h(p)}2 + \frac{1+(-1)^k}4
,
\end{gather*}
where $h(p)$ denotes the class number of the field
$\Q(\sqrt{-p})$. The occurrence of the class number is due to the sum
of the $\lambda_j\not=0$ in the general dimension formula, which, for
the binary $F$ in question, becomes a sum over $r_j/p$, where $0<
r_j<p$ runs through all quadratic residues modulo $p$.  For $p=3$ we
have to replace the third term on the right by
$0,0,-\frac13,\frac13,\frac13,-\frac13$ accordingly as $k$ modulo $6$
equals $0,1,2,3,4,5$, respectively, and $h(p)$
has to be replaced by $h(3):=\frac13$. It is not hard to show that $\dim
J_{1,F}=0$ (by Theorem~\ref{thm:singular-weight} and since $W(F)$
decomposes into two nontrivial irreducible characters of
dimension~$(p\pm1)/2$).

Table~\ref{tab:p_F} lists the {\it Hilbert-Poincar\'e polynomials}
$p_F$ for the first primes $p$; more precisely, it lists the
polynomials $\widetilde{p_F}=p_F-x^2\dim J_{2,F}$. For the dimension
of  the Jacobi forms of weight $2$ we have no clue to make a general
statement. It can be computed, however, for each $F$ as long as $p$ is
not too large.

Note that, for a fixed $k$, the dimensions $\dim J_{k,F}$ do only depend
on $\det(2F)$. This reflects the fact that, for binary $F$ with prime
discriminant, the isomorphism class of the determinant module $D_F$,
hence also the isomorphism class of $W(F)$, depend only on $\det(2F)$.
Indeed,~$D_F$ is
isomorphic to the quadratic module $\F_p$ equipped with the quadratic
form $Q(x)=x^2/p$. See section~\ref{Weil} for the terminology
used in this paragraph.

For deducing the given dimension formula
from Theorem~\ref{thm:dimension-formula}
it is useful to note that the $\SL$-action on $W(F)$ factors through
$\SL[\F_p]$ and $W(F)$ decomposes into
a $(p+1)/2$ and $(p-1)/2$ dimensional representation, which can be
easily identified by consulting a character table of $\SL[\F_p]$.

\begin{table}
\label{tab:p_F}
\begin{center}
\begin{tabular}[c]{|l|l|}
\hline
$p$&$\widetilde{p_F}(x)$\\
\hline
$3$&$x^9 + x^6 + x^4$\\
$7$&$x^{11} + x^{10} + x^9 + x^8 + x^7 + x^6 + x^4$\\
$11$&$x^{11} + x^{10} + 2 x^9 + 2 x^8 + x^7 + 2 x^6 + x^5 + x^4$\\
$19$&$x^{11} + 2 x^{10} + 3 x^9 + 3 x^8 + 3 x^7 + 3 x^6 + 2 x^5 + 2 x^4$\\
$23$&$x^{12} + 3 x^{11} + 3 x^{10} + 4 x^9 + 4 x^8 + 3 x^7 + 3 x^6 + x^5 + x^4$\\
$31$&$x^{12} + 3 x^{11} + 4 x^{10} + 5 x^9 + 5 x^8 + 5 x^7 + 4 x^6 + 2 x^5 + 2 x^4$\\
$43$&$2 x^{11} + 4 x^{10} + 6 x^9 + 7 x^8 + 7 x^7 + 7 x^6 + 5 x^5 + 4 x^4 + x^3$\\
$47$&$2 x^{12} + 5 x^{11} + 6 x^{10} + 8 x^9 + 8 x^8 + 7 x^7 + 6 x^6 + 3 x^5 + 2 x^4$\\
\hline
\end{tabular}\\[.5em]
\caption{
Hilbert-Poincar\'e polynomials
for binary $F$ with $\det(2F)=p$.}
\end{center}
\end{table}

For $p=3$ we find
$$
p_F(x)=\widetilde{p_F}(x) + c x^2 (1-x^4) (1-x^6)=c x^{12} + x^9 - c x^8 + (-c + 1) x^6 + x^4 + c x^2
,
$$
which implies $c=\dim J_{2,F}=0$. Accordingly, we have
$$
J_{*,F}
=
M_*\,\Psi_4
\oplus M_*\,\Psi_6
\oplus M_*\,\Psi_9
.
$$
for the (up to a constant factor) unique Jacobi forms
$\Psi_k\in J_{k,F}$ ($k=4,6,9$). To become more specific, choose
$2F=\begin{pmatrix}2&1\\1&2\end{pmatrix}$. Then one has:
$$
\Psi_9(\tau,z_1,z_2)=\vartheta(\tau,z_1)\,
\vartheta(\tau,z_1+z_2)\,\vartheta(\tau,z_2)\,\eta^{15}(\tau),
$$
where
$$
\vartheta
=
q^{\frac18}\big(\zeta^{\frac12}-\zeta^{-\frac12}\big)
\prod_{n\ge 1}
\big(1-q^n\big)
\big(1-q^n\zeta\big)
\big(1-q^n\zeta^{-1}\big),
\quad
\eta=q^{\frac1{14}}\prod_{n>1}\left(1-q^n\right).
$$
Here we use $q(\tau)=e^{2\pi i\tau}$ and $\zeta(\tau)=e^{2\pi
iz}$.
The formula for $\Psi_9$ follows immediately from the transformation
laws of the fundamental Jacobi form~$\vartheta$ (see the Appendix).
The nice product formula is due to the fact that, for fixed~$\tau$,
the theta function $\Psi_9(\tau,\cdot)$ has to be a multiple of
$\vartheta_{F,e}(\tau,\cdot)-\vartheta_{F,-e}(\tau,\cdot)$,
where $e\in\Z^2$, $e\not\in 2F\Z^2$, since $W(F)^c(+1)$ is one dimensional.
From this it is then clear that $\Psi_9(\tau,\cdot)$ vanishes
along $z_1=0$, $z_1+z_2=0$ and $z_2=0$, respectively, and hence
has to be divisible (in the ring of holomorphic functions)
by $\vartheta(\tau,z_1)\vartheta(\tau,z_1+z_2)\vartheta(\tau,z_2)$.

For $\Psi_4$ and $\Psi_6$ we find the closed formulas
$$
\Psi_k(\tau,z_1,z_2)
=
\sum_{n,a,b}
q^n\,e^{2\pi i(az_1+bz_2)}\,
\nu_{a^2-ab+b^2}\sum_{st=3n-a^2+ab-b^2}
s^{k-2}
\left[\left(\frac s3\right)-\left(\frac{t}3\right)\right]
.
$$
In the last formula formula we have $\nu_{N}=\frac12$ or $\nu_{N}=1$ accordingly as
$N$ is $1$ or $0$ modulo $3$. 
The outer sum
is over all integers $n,a,b$ such that
$$n-(a^2-ab+b^2)=n-F[(a,b)^t]/4\ge 0,$$
the inner sum is over all positive integers $s,t$ satisfying the
given identity,
and for
$a_k(N):=\sum_{st=N}
s^{k-2}
\left[\left(\frac
s3\right)-\left(\frac{t}3\right)\right]$ we use the convention
$a_k(0)=\frac12 L\left(2-k,\left(\frac \cdot3\right)\right)
$, in particular, $a_4(0)=-\frac19$ and $a_6(0)=\frac13$.

The
formulas for $\Psi_4$ and $\Psi_6$ can be deduced as follows. It is
well-known \cite[Thm.~1,~p.~333]{Kl 46} that $\vartheta_{F,0}$ lies in
$J_{1,F}\left(\Gamma_0(3),\left(\frac\cdot3\right)\right)$, hence
the sum $\vartheta_{F,0}+\vartheta_{F,e}+\vartheta_{F,-e}$
($e\in\Z^2$, $e\not \in 2F\Z^2$) defines an element of
$J_{1,F}\left(\Gamma^0(3),\left(\frac\cdot3\right)\right)$ (as follows
from $S^{-1}\Gamma_0(3)S=\Gamma^0(3)$, and from the action of the
$S$-matrix on the Weil representation $W(F)$ (see
section~\ref{Weil})). But then the function $h:=h_0+h_e+h_{-e}$ derived
from the decomposition
$\Psi_k=h_0\vartheta_{F,0}+h_e\vartheta_{F,e}+h_{-e}\vartheta_{F,-e}$
is a modular form of weight $k-1$ on $\Gamma^0(3)$ of nebentypus
$\left(\frac\cdot3\right)$ with a Fourier development of the form
$h=\sum_{n\equiv0,-1\bmod3} a(N)\,q^{N/3}$. For $k=4,6$ the
spaces of such modular forms are one dimensional, respectively; in fact,
these spaces are spanned by the Eisenstein series with Fourier
coefficients~$a_k(N)$. (Actually, for even $k$, the map
$\sum h_x\vartheta_{F,x}\mapsto \sum h_x$ defines an isomorphism of
$J_{k,F}$ and the space of modular forms of weight $k$ on $\Gamma^0(3)$
with real quadratic nebentypus and
with the aforementioned
special Fourier development. This will be proved elsewhere.
In particular, the formulas for $\Psi_k$
define an element of $J_{k,F}$ for all even $k$.)  For
calculating $\Psi_k$ from~$h$ note that $h_e=h_{-e}$ as follows from
the eveness of $\Psi_k(\tau,z)$ in $z$, as follows in turn from the
invariance of $\Psi_k$ under the modular transformation $-1$.

\section{Weil representations}
\label{Weil}

The metaplectic cover $\Mpp$ of $\SL$ acts on
$\Theta_F$, the space spanned by the $\det(2F)$ many theta series
$\vartheta_{F,x}$ (see the Appendix for the definition of the
metaplectic cover).  This is simply the algebraic restatement of the
well-known transformation formulas for theta series relating, for any
$A$ in the modular group, the series $\vartheta_{F,x}(A\tau,z)$ to the
series $\vartheta_{F,x}(\tau,z)$.  These transformation formulas stem
back to Jacobi and were treated by many authors. However,  a good
reference suited to our discussion is \cite{Kl 46}.  To focus
on the algebraic nature of the associated representation of $\Mpp$ we
use the notion of Weil representations associated to finite quadratic
modules.

Recall that a finite quadratic module $M$ is a a finite abelian group
$M$ endowed with a quadratic form $Q_M:M\rightarrow \Q/\Z$.
Here a
quadratic form is, by definition, a map such that $Q_M(ax)=a^2Q_M(x)$
for all $x$ in $M$ and all integers~$a$, and such that the application
$B_M(x,y):=Q_M(x+y)-Q_M(x)-Q_M(y)$ defines a $\Z$-bilinear map
$B_M:M\times M\rightarrow \Q/\Z$. All quadratic modules occurring in
the sequel will be assumed to be non degenerate, i.e.~we assume that
$B_M(x,y)=0$ for all $y$ is only possible for~$x=0$.

\begin{Definition}[Weil representation]
  The Weil representation $W(M)$ associated to a finite quadratic module~$M$
  is the $\Mpp$-module whose underlying vector space
  is $\C[M]$, the complex vector space of all formal linear
  combinations $\sum_x \lambda(x)\,\pt x$, where $\pt x$, for $x\in M$,
  is a symbol, where $\lambda(x)$ is a complex number and where the sum
  is over all $x$ in $M$, and where the action of $(T,w_T)$ and
  $(S,w_S)$ on $\C[M]$ is given by the formulas
  \begin{equation*}
	 \label{eq:action}
	 \begin{split}
		(T,w_T)\,\pt x &= e(Q_M(x))\,\pt x\\
		(S,w_S)\,\pt x &= \sigma\,|M|^{-\frac12}\sum_{y\in M} \pt y\, e(-B_M(y,x))
		,
	 \end{split}
  \end{equation*}
  respectively, where
  $$
  \sigma=\sigma(M)=|M|^{-\frac12}\sum_{x\in M} e(-Q_M(x))
  .
  $$
\end{Definition}

That the action of the generators
$$
(T,w_T)=\left(\begin{pmatrix}1&1\\0&1\end{pmatrix},1\right),
\quad
(S,w_S)=\left(\begin{pmatrix}0&-1\\1&0\end{pmatrix},\sqrt \tau\right)
$$
of $\Mpp$ can indeed
be extended to an action of $\Mpp$ on $\C[M]$ is discussed in more detail
in \cite{Sko 07}. For the moment we do not know any reference
for this fact and shall explain its (simple) proof elsewhere~\cite{Sko 08}.
The action of
$\Mpp$ factors through $\SL$ if and only
if~$\sigma^4=1$~\cite{Sko 08};
in general, $\sigma$ is an eighth root of unity. That the
actions of $(S,w_S)$ and $(T,w_T)$ (or rather $S$ and $T$)
define an action of
$\SL$ if $\sigma^4=1$ is well-known
(cf.~e.g.~\cite{N 76}).

The standard example for a quadratic module is the determinant
module~$D_F$ of a symmetric non degenerate half integral matrix $F$.
The quadratic module has $D_F=\Z^n/2F\Z^n$
as underlying abelian group, and the quadratic form on $D_F$ is the one
induced by the quadratic form $\Z^n \ni x\mapsto \frac 14 F^{-1}[x]$.

The $\Mpp$-module $W(F):=W(D_F)$ is the one occurring in the formulas
of the first section.  The associated right module $W(F)^*$ is isomorphic
(as $\Mpp$-module) to $\Theta_F$. An isomorphism is given by the map
$$
W(F)^*\rightarrow\Theta_F,
\quad
\lambda\mapsto\sum_{x\in\Z^n/2F\Z^n} \lambda(\pt x)\,\vartheta_{F,x}
.
$$
That this map is a $\Mpp$-right module isomorphism is clear from
from the formulas for the action of $(S,w_S)$ and $(T,w_T)$ on $W(F)$
as given above
and on $\Theta_F$ as given in~\cite[(1.12), p.~320]{Kl 46}.

In our context, the most important notion related to a quadratic module
is its associated space of invariants.
\begin{Definition}[Invariants associated to a quadratic module]
  For a quadratic module $M$, we use
  $\Inv(M)$ for
  the subspace of elements in $W(M)$ which are invariant under the action of $\Mpp$
\end{Definition}

We already saw in section~\ref{Jacobi} that $J_{\frac n2,F}(\SL)$, for
index $F$ of even size~$n$, is isomorphic to $\Inv(D_F)$ (cf.~the remark after
Theorem~\ref{thm:singular-weight}) and that a similar statement holds true
for the critical weight case $J_{\frac {n+1}2,F}(\SL)$ if $n$ is odd.
In fact, it can be shown that, for
any $\Gamma$ and character $\chi$ which is trivial on a congruence subgroup, 
the spaces $J_{\frac n2,F}(\Gamma,\chi)$ or $J_{\frac {n+1}2,F}(\Gamma,\chi)$,
accordingly as $n$ is even or odd, are isomorphic to certain natural subspaces
of invariants of quadratic modules \cite{Sko 07}. Thus, the clue to Jacobi forms
of critical weight (and weight $\frac n2$ for even $n$) is the study
of invariants of Weil representations. A rather involved illustration
for this statement can be found in \cite{Sko 07}, where the spaces
$J_{1,m}(\SL,\varepsilon^8)$
are explicitly determined for all~$m$
 (see the Appendix for the character $\varepsilon$).

The question for the invariants of a quadratic module seems to be
subtle.  We do not even know a reasonable (say, easily computable)
criterion to answer in general the
question which quadratic modules possess invariants and which
not. However, it seems that the nature of invariants depend on the
Witt class of the underlying quadratic module.

For stating this more precisely note that, for an isotropic submodule
$N$ of a quadratic module $M$\footnote{%
$N$ is called {\it isotropic
submodule of} $M$ if $Q_M(x)=0$ for all $x\in N$.
},
the dual module
$N^\perp=\{y\in M:B_Q(N,y)=0\}$ contains $N$ and $Q_M$ induces a
quadratic form ${\overline Q_M}$ on $N^\perp/N$. In fact, $N^\perp/N$,
equipped with the quadratic form induced by $Q_M$ becomes a
(non degenerate) quadratic module. We call $M$ Witt-zero\footnote{%
Some authors use the terminology {\it weakly metabolic} instead.
}
if $M$ contains
an isotropic {\it self-dual} module $N$ i.e.~an isotropic module $N$
such that $N=N^\perp$.

Two quadratic modules $M_1$ and $M_2$ are called
Witt-equivalent if they contain isotropic submodules $N_1$ and $N_2$,
respectively, such that $N_1^\perp/N_1$ and $N_2^\perp/N_2$ are
isomorphic as quadratic modules. It is not hard to check that this
defines indeed an equivalence relation. The set of equivalence
classes~$[M]$ of all (non degenerate finite) quadratic modules $M$ becomes a
group via the operation $[M]+[N]:=[M\perp N]$, where $M\perp N$ is the quadratic
module with underlying group $M\oplus N$ and quadratic form $Q(x\oplus
y)=Q_M(x)+ Q_N(y)$. Adopting the notation in \cite{Sch 84} we denote
this group by $WQ$. Similarly, for any prime $p$, we can define the
{\it Witt group} $WQ(p)$ of all (finite non degenerate) quadratic
modules whose order is a power of $p$.  If $p$ is odd, the group
$WQ(p)$ is of order $4$, whereas $WQ(2)$ is of order 16 \cite{Sko 08},
\cite{Sch 84} (This is not hard to prove: If $N$ is a maximal
isotropic submodule of $M$ then $P:=N^\perp/N$ is anisotropic,
i.e. satisfies $Q_P(x)\not=0$ for all $x\not=0$, and $M$ and~$P$
span the same Witt class; thus it suffices to count the isomorphism
classes of anisotropic quadratic modules which are $p$-groups).

Every quadratic module $M$ possesses the {\it primary decomposition}
$$
M=\perp_p M(p)
.
$$
Here $p$ runs through all primes and $M(p)$ denotes the $p$-part of
$M$, i.e. the submodule of all $x$ in $M$ whose order is a power of
$p$, equipped with the restriction of $Q_M$. The primary decomposition
gives rise to a canonical isomorphism
$$
WQ
\cong 
\coprod_p WQ(p)
.
$$
Moreover, via the functorial identity $W(M\perp N)=W(M)\otimes W(N)$
it implies the following lemma.
\begin{Lemma}
  $$
  \Inv(M)
  \cong
  \bigotimes_{p|l} \Inv\big(M(p)\big)
  .
  $$
\end{Lemma}

Thus, to study Weil invariants of quadratic modules we can restrict to
quadratic modules which are $p$-groups. For these we have
the following theorem \cite{N-R-S 06}.
\begin{Theorem}[N-R-S 06]
Let $M$ be a quadratic of prime power order whose Witt class vanishes.
Then $\Inv(M)$ is different from
zero. Moreover, $\Inv(M)$ is generated by all $I_U=\sum_{x\in U}\pt x$,
where $U$ runs through the isotropic self-dual subgroups of $M$.
\end{Theorem}

However, for most of the cases which we are interested in the
assumption of the last theorem fails for some prime $p$. Indeed, for a
positive half integral $F$, the module $D_F$ is Witt-zero if and only
if $\det(2F)$ is a square and $\sym{rank}(F)\equiv 0\bmod 8$ (this is
e.g.~explaind in \cite{Sko 07} and follows from Milgram's formula and
the theory of $\sigma$-invariants explained in loc.cit.).  But if
$D_F$ is not Witt-zero, then one of its primary constituent must not be
Witt-zero either and we cannot apply the last theorem.

In \cite{Sko 07} the interested reader finds some more results on the
Weil invariants of quadratic modules. However, we do not yet have a
complete theory. On the other hand it seems to us that the
characterization of the Weil invariants of quadratic modules and the
search for a dimension formula for the space spanned by the Weil
invariants is a a major problem in the context of Jacobi forms of
small weight and their application to the theory of Siegel and
orthogonal modular forms.

\section*{Appendix}
\label{sec:Appendix}

For the reader's convenience we recall in this section those technical
details concerning Jacobi forms which we tacitly used in this article.
  
The metaplectic double cover $\Mpp$ of $\SL$, consists of all pairs
$(A,\pm w_A)$, where $A=\begin{pmatrix}a&b\\c&d\end{pmatrix}\in\SL$, and
where $w_A$ is the holomorphic function on $\HP$ such that
$w_A(\tau)^2=c\tau+d$ and such that
$-\pi/2<\sym{Arg}(w_A(\tau))\le \pi/2$. The composition
law is given by $(A,w)\cdot(B,v)=\big(AB,w(B\tau)v(\tau)\big)$.
If~$F$ is a real
symmetric, half integral  $n\times n$~matrix and $k$ a half integral
integer we have the action of $\Mpp$ on functions $\phi$ defined on
$\HP\times\C^n$ given by the formula:
$$
\phi|_{k,F}(A,w)(\tau,z)
=
\phi\left(A\tau,\frac z{c\tau+d}\right)\,
w(\tau)^{-2k}\,
\e\left(\frac {-c F[z]}{c\tau+d}\right).
$$
If $k$ is integral this action factors through an action of $\SL$ and
we write $\phi|_{k,F}A$ for $\phi|_{k,F}(A,w)$.  Furthermore, if $F$
is half integral, we have the action of $\Z^n\times \Z^n$ on functions on
$\HP\times\C^n$ defined by
$$
\phi|_{k,F}[\lambda,\mu](\tau,z)
=
\phi(\tau,z+\lambda\tau +\mu)\,
\e\left(\tau F[\lambda] + 2z^tF\lambda\right)
$$

\begin{Definition}(Jacobi forms)
  Let $F$ be a symmetric, half integral positive definite $n\times
  n$~matrix, let $k$ be an integer and let $\Gamma$ be a subgroup of
  finite index in $\SL$ and $\chi$ be a linear character of $\Gamma$. A
  {\sl Jacobi form of weight~$k$ and index $F$ on $\Gamma$ with
  character~$\chi$} is a holomorphic function $\phi:\HP\times
  \C^n\rightarrow \C$ such that the following holds true:
  \begin{enumerate}
  \item[(i)]
    For all $A \in \Gamma$ one has $\phi|_{k,F}A = \chi(A)\phi$,
    and for all $g\in \Z^n\times\Z^n$ one has $\phi|_{k,F}g = \phi$.
  \item[(ii)]
    For all $A\in\SL$ the function $\phi|_{k,F}A$
    possesses a Fourier expansion of the form
    $$
    \phi|_{k,F}A = \sum_{\begin{subarray}{c}l\in\Q,\; r\in\Z^{n}\\
		  4l-F^{-1}[r]\ge 0\end{subarray}} c(l,r)\,q^l\,\e(z^tr).
    $$
  \end{enumerate}
The space of these Jacobi forms is denoted by $J_{k,F}(\Gamma,\chi)$.
\end{Definition}

Similarly, we can define Jacobi forms of half integral weight $k$ on
subgroups~$\Gamma$ of $\Mpp$. As pointed out earlier there should also
be a theory of Jacobi forms of index $F$ where the diagonal elements of
$F$ are half integral. However, then the transformation law with
respect to $\Z^n\times\Z^n$ would involve a character of $\Z^n\times\Z^n$
of order $2$. Moreover, the theory of these forms is, in a sense, included in
the theory of Jacobi forms of half integral index by the map
$\phi(\tau,z)\mapsto \phi(\tau,2z)$, which maps Jacobi forms of index $F$
(not necessarily with integral diagonal) to Jacobi forms of index $2F$.

The simplest nontrivial Jacobi form is the function
$$
\vartheta(\tau,z)
=
\sum_{n\in\Z}\left(\frac{-4}n\right)q^{\frac{n^2}8}\zeta^{\frac n2}
=
q^{\frac18}\big(\zeta^{\frac12}-\zeta^{-\frac12}\big)
\prod_{n\ge 1}
\big(1-q^n\big)
\big(1-q^n\zeta\big)
\big(1-q^n\zeta^{-1}\big)
$$
(the second identity follows from the Jacobi triple product identity).
In fact, for $(A,w)$ in $\Mpp$ and integers $\lambda,\mu$, this
function satisfies the following transformation formulas
\cite[p.~145]{Sko 92}:
\begin{gather*}
\vartheta|_{\frac12,\frac12}(A,w)=\epsilon^3(A,w)\vartheta\\
\vartheta|_{\frac12,\frac12}(\lambda,\mu)=(-1)^{\lambda+\mu}\vartheta
,
\end{gather*}
where $\epsilon(A,w)=\eta(A\tau)/w(\tau)\eta(\tau)$. This is not a
Jacobi form in the strict sense as defined above, but one of {\it
index $\frac 12$}. In any case, $\vartheta(\tau,2z)$ is an element of
$J_{\frac12,2}(\Mpp,\epsilon^3)$.

Finally, we remark, that a Jacobi form of index $F$ defines indeed a
holomorphic section of a positive line bundle on $\C/\Z^n\tau+\Z^n$
with the Chern class associated to $F$ as
explained in section~\ref{Jacobi}. A hermitian form on this line bundle
is induced by the function $H(\tau,z)=e^{-4\pi F[\Im(z)]/\Im(\tau)}$. 
In fact, if $\phi$ and $\psi$ are Jacobi forms
of index $F$, then $\phi H\overline{\psi}$ is invariant under
$z\mapsto z+\lambda\tau+\mu$.  If we write $F=\left(f_{p,q}\right)$,
the Chern class of the line bundle in question is thus given by
$$
\frac1{2\pi i}\partial \overline\partial \log H
=
\frac 1{i\Im(\tau)}\sum_{p,q}f_{p,q}\,dz_p\wedge d\overline z_q
.
$$


\bigskip
\begin{flushleft}
Nils-Peter Skoruppa\\
Fachbereich Mathematik,
Universit\"at Siegen\\
Walter-Flex-Stra{\ss}e 3,
57068 Siegen, Germany\\
http:/\hskip-3pt/www.countnumber.de\\
\end{flushleft}

\end{document}